\documentclass{ifacconf}

\usepackage[caption=false]{subfig}
\usepackage[dvips]{graphicx}
\graphicspath{{figures/}}
\usepackage{amsmath, amssymb}
\usepackage{tikz} \usetikzlibrary{calc}
\usepackage{psfrag}
\usepackage{pgfplots}
\usepackage{algorithm}     
\usepackage{algpseudocode} 
\usepackage{arydshln}

\DeclareMathOperator*{\minimize}{minimize}

\DeclareMathOperator*{\subject}{subj.\ to}

\DeclareMathOperator*{\argmin}{arg\ min}

\DeclareMathOperator*{\diag}{diag}

\DeclareMathOperator*{\rank}{rank}

\newcommand{\nt}{\text{nt}}

\definecolor{red}{rgb}{1,0,0}

\begin{document}

\begin{frontmatter}

\title{Distributed Interior-point Method for Loosely Coupled Problems\thanksref{footnoteinfo}} 
\vspace*{-7pt}
\thanks[footnoteinfo]{This work has been supported by the Swedish Department of Education within the ELLIIT project.}
\author[First]{Sina Khoshfetrat Pakazad} 
\author[First]{Anders Hansson}
\author[Second]{Martin S. Andersen} 
%
\address[First]{Division of Automatic Control, Department of Electrical Engineering, Link\"oping University, Sweden. Email: \{sina.kh.pa, hansson\}@isy.liu.se.} 
\address[Second]{Department of Applied Mathematics and Computer Science, Technical University of Denmark, Denmark. (e-mail: mskan@dtu.dk)}                                         


\begin{abstract}                          
In this paper, we put forth distributed algorithms for solving loosely coupled unconstrained and constrained optimization problems. Such problems are usually solved using algorithms that are based on a combination of decomposition and first order methods. These algorithms are commonly very slow and require many iterations to converge. In order to alleviate this issue, we propose algorithms that combine the Newton and interior-point methods with proximal splitting methods for solving such problems. Particularly, the algorithm for solving unconstrained loosely coupled problems, is based on Newton's method and utilizes proximal splitting to distribute the computations for calculating the Newton step at each iteration. A combination of this algorithm and the interior-point method is then used to introduce a distributed algorithm for solving constrained loosely coupled problems. We also provide guidelines on how to implement the proposed methods efficiently and briefly discuss the properties of the resulting solutions. 
\end{abstract}

\end{frontmatter}

\section{Introduction}
\vspace*{-7pt}
Distributed optimization methods are specially crucial in the absence of a centralized computational unit or when this unit lacks the necessary computational power for solving the problem at all or in a timely manner. In either case, distributed algorithms enable us to solve the problem without the need for a central computational unit and through collaboration of a number of computational agents. Such methods also allow us to solve the problem by solving subproblems that are easier and possibly faster to solve. Distributed methods for solving optimization problems have been studied for many years and there are different approaches for devising such algorithms, see e.g., \cite{ber:97,eck:89,boyd:11,ned:09,ned:10}.

One of the most commonly used approaches for devising distributed optimization algorithms is based on applying gradient/subgradient methods directly to the problem, see e.g., \cite{ned:09,ned:10}. The resulting algorithms from such approaches are usually simple and easy to implement, however, they are very sensitive to the scaling of the problem and suffer from very slow convergence, \cite{ber:97}. Another approach for constructing distributed algorithms is based on decomposition techniques and proximal splitting methods. To this end, the optimization problem is first decomposed and then a splitting method of choice, e.g., the alternating direction method of multipliers (ADMM), is applied to the decomposed problem. This will define the computational routines that each agent should perform locally and also will set the communication/collaboration protocol among the agents, e.g., see \cite{ber:97,eck:89,boyd:11,com:11}. Even though the resulting algorithms using this approach tend to be more complicated than the previous type of algorithms, they enjoy a faster rate of convergence for certain classes of optimization problems, e.g., when the objective function of the equivalent unconstrained reformulation of the problem, has two terms and/or is strongly convex, \cite{gol:12,gols:12}. For more general classes of problems, however, for instance when the objective function has more than two terms of which several are non-smooth, e.g., indicator functions, they can perform very poorly or may even fail to converge, e.g., see \cite{che:13}. Although there exist modifications to splitting methods that allow us to apply them to more general classes of optimization problems, the resulting algorithms can become overly complicated to implement (particularly distributedly), \cite{gol:12,han:12,hon:12}, or the local computations can still be considerable, e.g., they may require each agent to solve an equality/inequality constrained problem at each iteration, \cite{sum:12,ohl:13}.

In order to allay the aforementioned issues, there has been a recent interest in studying the possibility of using second order methods for designing distributed optimization algorithms, e.g., see \cite{Chu:11,wei:13,nec:09}. For instance, in \cite{wei:13}, the authors propose a distributed Newton method for solving a network utility maximization problem. The cost function for such problems comprise of a summation of several terms where each term depends on a single scalar variable. This structure allows the authors to employ a matrix splitting method which in turn enables them to compute the Newton directions distributedly. However, this method relies on the special structure in the considered problem and hence can only be used for problems with such structure. In \cite{nec:09} the authors propose a distributed optimization method based on an interior-point method. The introduced algorithm is obtained by first performing a Lagrangian decomposition of the problem and then solving the subproblems using interior-point methods, efficiently. However, in the proposed algorithm, the computational cost for solving the subproblems can still be considerable. The authors in \cite{Chu:11} propose a distributed Newton method for solving coupled unconstrained quadratic problems, which is used for anomaly detection in large populations. This distributed method is only applicable to unconstrained quadratic problems.  

In this paper, we investigate the possibility of utilizing Newton's and interior-point methods for designing distributed  algorithms for solving loosely coupled optimization problems. Notice that this type of problems constitute a more general class of problems than the ones in the above mentioned papers. For this purpose, we first exploit the coupling in the problem using consistency constraints and use proximal splitting methods to compute the Newton directions in a distributed manner. This then enables us to propose distributed implementations of the Newton method, which can be used for solving unconstrained loosely coupled problems. In order to solve constrained loosely coupled problems, we also put forth a distributed interior-point method which relies on the proposed distributed Newton method. Furthermore, in contrast to the methods proposed in \cite{ned:09}, in all the proposed algorithms each of the agents is only required to solve an unconstrained or equality constrained quadratic program at each iteration.

\vspace*{-5pt}
\subsection*{Notation}
\vspace*{-5pt}
We denote by $\mathbb R$ the set of real scalars and by $\mathbb R^{n\times m}$ the set of real $n\times m$ matrices. The transpose of a matrix $A$ is denoted by $A^T$ and the column and null space of this matrix is denoted by $\mathcal{C}(A)$ and $\mathcal N(A)$, respectively. We denote the set of positive integers
$\{1,2,\ldots,p\}$ with $\mathbb{N}_p$. Given a set $J \subset \{1,2,\ldots,n\}$, the matrix $E_J \in \mathbb{R}^{|J|\times n}$ is the $0$-$1$ matrix that is obtained by deleting the rows indexed by $\mathbb{N}_n \setminus J$ from an identity matrix of order $n$, where $|J|$ denotes the number of elements in set $J$. This means that $E_Jx$ is a $|J|$- dimensional vector with the components of $x$ that correspond to the elements in $J$, and we denote this vector with $x_J$. With $x^{i,(k)}_l$ we denote the $l$th element of vector $x^i$ at the $k$th iteration. Also given vectors $x^i$ for $i= 1, \dots, N$, the column vector $(x^1, \dots, x^N)$ is all of the given vectors stacked. The orthogonal projection onto a set $\mathcal C$ is denoted by $P_{\mathcal C}$.
\vspace*{-5pt}
\section{Newton's Method for Solving Equality Constraint Problems}\label{sec:New}
\vspace*{-5pt}
Consider the following equality-constrained optimization problem 
\begin{subequations}\label{eq:EOP}
\begin{align}
\minimize_x    & \quad F(x)\\
\subject   &\quad A x = b,
\end{align} 
\end{subequations}
where $F: \mathbb R^n \rightarrow \mathbb R$ is convex and twice differentiable, and $A \in \mathbb R^{p\times n}$ with $\rank(A) = p < n$. This problem can be solved iteratively using Newton's method, \cite[Ch. 10]{boyd:04}, where at each iteration the variables are updated as $x^{(l+1)} = x^{(l)} + \alpha^{(l)}\Delta x^{(l)}_{\nt}$, with $\Delta x^{(l)}_{\nt}$ and $\alpha^{(l)}$ denoting the so-called Newton direction and its corresponding step size, respectively. The Newton direction $\Delta x^{(l)}_{\nt}$ is defined as the solution of  
\begin{subequations}\label{eq:QApp}
\begin{align}
\minimize_{\Delta x} & \quad F(x^{(l)}) + \nabla F(x^{(l)})^T \Delta x + \frac{1}{2}  \Delta x^T \nabla^2F(x^{(l)}) \Delta x\\
\subject & \quad \quad A(x^{(l)} + \Delta x) = b, \label{eq:QApp_b}
\end{align}
\end{subequations}
which constitutes a quadratic approximation of \eqref{eq:EOP} at $x^{(l)}$, \cite[Ch. 10]{boyd:04}. Notice that if we assume that $Ax^{(l)} = b$ for all $l\geq 0$, then the constraint in~\eqref{eq:QApp_b} would be $A\Delta x = 0$. The corresponding step size $\alpha^{(l)}$ is commonly computed using either exact or backtracking line search methods. The stopping/termination criterion for this iterative algorithm is based on the so-called Newton decrement which is defined as
\begin{align}\label{eq:ND}
\lambda^{(l)} = \sqrt{(\Delta x^{(l)}_{\nt})^T \nabla^2F(x^{(l)})\Delta x^{(l)}_{\nt}},
\end{align}
and it provides an estimate of the sub-optimality of the iterates. Hence, Newton's method terminates when the Newton decrement falls below a given threshold $\epsilon$. This iterative scheme is summarized in Algorithm \ref{alg:alg1}.
\begin{algorithm}[tb]
\caption{Newton's Method}\label{alg:alg1}
\begin{algorithmic}[1]
\small
\State{Given $l =0$, $\epsilon_{\nt}>0$ and $x^{(0)}$ such that $Ax^{(0)} = b$}
\Repeat
\State{Compute $\Delta x^{(l)}_{\nt}$ by solving \eqref{eq:QApp}}
\State{Compute $\lambda^{(l)}$}
\If  {$(\lambda^{(l)})^2/2 \leq \epsilon_{\nt}$} 
\State{$x^\ast = x^{(l)}$} 
\State {Terminate the iterations} 
\EndIf 
\State {Compute the step size $\alpha^{(l)}$ using line search}
\State  $x^{(l+1)} =  x^{(l)} + \alpha^{(l)}\Delta x^{(l)}_{\nt}$
\State{$l = l+1$}
\Until{iterations are terminated}
\normalsize
\end{algorithmic}
\end{algorithm}
\vspace*{-5pt}
\section{Interior-point Method}\label{sec:In}
\vspace*{-5pt}
The optimization problem in \eqref{eq:EOP} can be extended to include inequality constraints as
\begin{subequations}\label{eq:EOPInC}
\begin{align}
\minimize_x   &\quad F(x)\\
\subject   &\quad A x = b\\
& \quad g_i(x) \leq 0 \quad i = 1, \dots, m,\label{eq:EOPInC_c}
\end{align}
\end{subequations}
where $g_i:\mathbb R^n \rightarrow \mathbb R$ are all convex and twice differentiable. Also assume that the problem admits a strictly feasible solution, i.e., there exists $ x $ such that $Ax = b$ and $g_i(x) < 0$ for all $i = 1, \dots, m$. This allows us to solve \eqref{eq:EOPInC} using an interior-point method which requires solving a sequence of equality-constraint problems, defined as 
\begin{subequations}\label{eq:EOPInB}
\begin{align}
\minimize_x    & \quad tF(x) + \sum_{i=1}^{m} \phi_i(x)\\
\subject   &\quad A x = b,
\end{align} 
\end{subequations}
for increasing values of $t>0$. The functions $\phi_i(x) := -\log(-g_i(x))$ are so-called logarithmic barrier functions which (since $\phi(x) \rightarrow \infty$ as $g_i(x) \rightarrow 0$) virtually create barriers that prevent $x$ from violating any of the inequality constraints in \eqref{eq:EOPInC_c}. Notice that for every fixed value of $t$ the problem in \eqref{eq:EOPInB} is precisely in the form of \eqref{eq:EOP}, and it can be solved using Newton's method, as is described in Algorithm \ref{alg:alg1}. After the Newton method has converged to a solution for a given $t$, $x^\ast(t)$, then within the interior-point method procedure, we gradually increase $t$ and perform the same procedure again. This is done until $m/t$ is below a certain threshold. This procedure is summarized in Algorithm~\ref{alg:alg2}.
\begin{algorithm}[tb]
\caption{Interior-point Method}\label{alg:alg2}
\begin{algorithmic}[1]
\small
\State{Given $q = 0$, $t^{(0)} >0$, $\mu>1$, $\epsilon_{p}>0$ and feasible $x^{(0)}$}
\Repeat
\State{Compute $x^\ast(t^{q})$ by solving \eqref{eq:EOPInB} using Alg. \ref{alg:alg1} starting at $x^{(q)}$}
\State{$x^{(q+1)} = x^\ast(t^{q})$}
\If  {$m/t^{(q)} < \epsilon_p$} 
\State{$x^\ast = x^{(q+1)}$} 
\State {Terminate the iterations} 
\EndIf 
\State{$t^{(q+1)} = \mu t^{(q)}$}
\State{$q = q+1$}
\Until{iterations are terminated}
\normalsize
\end{algorithmic}
\end{algorithm}
Next, we define loosely coupled problems and describe how we can take advantage of distributed computations  for solving such problems using algorithms~\ref{alg:alg1} and \ref{alg:alg2}. 
\vspace*{-10pt}
\section{ Unconstrained Loosely Coupled Optimization Problems}
\vspace*{-8pt}
\subsection{A Definition}\label{sec:Definition}
\vspace*{-5pt}
Consider the following unconstrained optimization problem 
\begin{align}\label{eq:UnconstGen}
\minimize \quad F_1(x) + \dots + F_N(x),
\end{align}
where the convex and twice differentiable functions $F_i$ for $i = 1, \dots, N$, only depend on a small subset of the elements of the variable $x$. Particularly, let us denote the ordered set of indices of variables that appear in $F_i$ by $J_i$. We also denote  the ordered set of indices of terms in the cost function that depend on $x_i$ by $\mathcal I_i$, i.e., $\mathcal{I}_i = \{k \ | \ i \in J_k\}$. We call an optimization problem loosely coupled if $|\mathcal I_i|\ll N$ for all $i = 1, \dots, n$. This problem can be rewritten as the following equality-constrained optimization problem 
\begin{subequations}\label{eq:DEOP}
\begin{align}
\minimize_{S,x}   \quad & f_1(s^1) + \dots + f_N(s^N)\\
\subject \quad  & \bar E x = S,\label{eq:DEOPb}
\end{align}
\end{subequations}
where $S = (s^1, \dots, s^N)$ and $\bar{E}  = \begin{bmatrix} E_{J_1}^T   &\cdots &  E_{J_N}^T \end{bmatrix}^T$, with $E_{J_i}$ a $0$-$1$ matrix that is obtained from an identity matrix of order $n$ by deleting the rows indexed by $\mathbb{N}_n \setminus J$. We refer to the constraints in \eqref{eq:DEOPb} as consistency constraints. Also $f_i$s are lower dimensional descriptions of $F_i$s such that $F_i(x) = f_i(E_{J_i}x)$ for all $x \in \mathbb R^n$ and $i = 1, \dots, N$. Notice that if we define $F(S,x) =  f_1(s^1) + \dots + f_N(s^N)$ and $A = \begin{bmatrix} -I & \bar E \end{bmatrix}$, the problem in~\eqref{eq:DEOP} is in the same form as the problem in \eqref{eq:EOP}.  We can now form the quadratic approximation in~\eqref{eq:QApp} for this problem, as below
\small
\begin{subequations}\label{eq:QAppLoose}
\begin{align}
\minimize_{\Delta s^1, \dots \Delta s^1, \Delta x} & \quad \sum_{i = 1}^N f_i(s^{i,(l)}) + \nabla f_i(s^{i,(l)})^T \Delta s^i + \notag\\ &\hspace{30mm} \frac{1}{2}  (\Delta s^i)^T \nabla^2f_i(s^{i,(l)}) \Delta s^i\\
\subject & \quad \bar E(x^{(l)} + \Delta x) = S^{(l)} + \Delta S,\label{eq:QAppLooseb}
\end{align}
\end{subequations}
\normalsize

where $\Delta S = (\Delta s^1, \Delta s^2, \dots, \Delta s^N)$. Assuming that $S^{(l)} \in \mathcal C(\bar E)$, \eqref{eq:QAppLooseb} only requires $\Delta S \in \mathcal C(\bar E)$. This problem can then be rewritten as the following unconstrained non-smooth optimization problem 
\small
\begin{multline}\label{eq:QAppLooseNon}
\minimize_{\Delta s^1, \dots \Delta s^1}  \quad \sum_{i = 1}^N f_i(s^{i,(l)}) + \nabla f_i(s^{i,(l)})^T \Delta s^i +\\ \frac{1}{2}  (\Delta s^i)^T \nabla^2f_i(s^{i,(l)}) \Delta s^i + \mathcal I_{\mathcal C} (\Delta S),
\end{multline}
\normalsize

where $\mathcal I_{\mathcal C}$ is the indicator function for the column space of $\bar E$. Considering the imposed structure in the cost function of the problem in \eqref{eq:QAppLooseNon}, it can be solved distributedly using so-called proximal splitting, \cite{com:11,boyd:11,eck:89}, which is the subject of the next section.
\vspace*{-5pt}
\subsection{Proximal Splitting Methods}
\vspace*{-5pt}
In many applications we are faced with convex problems of the form 
\begin{align}\label{eq:prox}
\minimize \quad T_1(x) + T_2(x)
\end{align}
where minimizing the joint problem is much harder than solving problems including only individual terms of the cost function, \cite{com:11,boyd:11,kho:13-1}. Proximal splitting algorithms enable us to solve the problem in~\eqref{eq:prox} by solving minimization problems that are based on individual terms in the cost function. This is done through the use of the so-called proximity operators of these terms which are defined as follows. Given a closed convex function $T$, the proximity operator for this function, $\text{prox}_T:\mathbb R^n \rightarrow \mathbb R$, is defined as the unique minimizer of 
\begin{align*}
\minimize_y \quad T(y) + \frac{1}{2} \| x-y \|^2.
\end{align*}
There are different classes of proximal splitting methods, \cite{eck:89}. In this paper we, however, only consider the Alternating Direction Method of Multipliers (ADMM) which has been extensively used recently for design of distributed algorithms in many applications, e.g., see \cite{boyd:11,ber:97}. A general description of ADMM for the optimization problem in \eqref{eq:prox} is given in Algorithm~\ref{alg:alg3}.
\begin{algorithm}[tb]
\caption{ADMM}\label{alg:alg3}
\begin{algorithmic}[1]
\small
\State{Given $k= 0$, $\rho > 0$, $\epsilon_{\text{pri}},\epsilon_{\text{dual}}>0$, $y^{(0)}$ and $v^{(0)}= 0$}
\Repeat
\State  $x^{(k+1)} = \text{prox}_{\frac{1}{\rho}T_1}(y^{(k)} + v^{(k)})$
\State  $y^{(k+1)} = \text{prox}_{\frac{1}{\rho}T_2}(x^{(k+1)} - v^{(k)})$
\State  $v^{(k+1)} =  v^{(k)} + (y^{(k+1)} - x^{(k+1)})$
\If  {$\| y^{(k+1)} - x^{(k+1)}\|^2 < \epsilon_{\text{pri}}$ \&\& $\| y^{(k+1)} - y^{(k)}\|^2 < \epsilon_{\text{dual}}$} 
\State {Terminate the algorithm} 
\EndIf 
\State{$k = k+1$}
\Until{algorithm is terminated}
\normalsize
\end{algorithmic}
\end{algorithm}
This algorithm can also accommodate a varying penalty parameter, $\rho$. This has shown to improve the convergence properties of the algorithm, and depending on the specifications of the problem there are different approaches for updating this parameter at each iteration, \cite{boyd:11,gha:13}. Also there exist several alternatives for the termination of the algorithm (the 6th step of Algorithm \ref{alg:alg3}), see e.g., \cite{boyd:11,ber:97}.  Proximal splitting algorithms and particularly ADMM can perform very well when applied to quadratic problems, \cite{bol:12}, and as we will see later they can be implemented very efficiently. Next, we describe how ADMM enables us to distribute the computations for solving \eqref{eq:DEOP}.
 \vspace*{-7pt} 
\subsection{Distributed Computation of the Newton Direction}\label{sec:DNew}
\vspace*{-7pt}
As was discussed in Section \ref{sec:Definition}, the problems in \eqref{eq:EOP} and~\eqref{eq:DEOP} have the same structure, and hence we can solve~\eqref{eq:DEOP} using the Newton method as described in Algorithm \ref{alg:alg1}. Recall that in the Newton method, in order to compute the Newton direction at each iteration, we need to solve the corresponding quadratic approximation of the problem. As was also shown in Section \ref{sec:Definition}, the quadratic approximation for the problem in \eqref{eq:DEOP}, can be equivalently rewritten as in~\eqref{eq:QAppLooseNon}. Notice that \eqref{eq:QAppLooseNon} is in the same format as \eqref{eq:prox}, where $T_1$ involves a summation of $N$ decoupled terms and $T_2$ corresponds to $\mathcal I_{\mathcal C}$. With these definitions for $T_1$ and $T_2$ we now apply ADMM to \eqref{eq:QAppLooseNon}. The 3rd and 4th steps of Algorithm \ref{alg:alg3} can then be written as 

\vspace*{-15pt}
\small
\begin{subequations}\label{eq:ADMM}
\begin{align}
\Delta S^{(k+1)} &= \text{prox}_{\frac{1}{\rho}T_1} (\Delta Y^{(k)}+v^{(k)})\notag\\
& = \argmin_{\Delta s^1, \dots, \Delta s^N} \quad \sum_{i = 1}^N \bar f_i(\Delta s^{i}) + \frac{\rho}{2}\| \Delta y^{i,(k)}+v^{i,(k)} - \Delta s^{i}  \|^2\label{eq:ADMMa} \\
\Delta Y^{(k+1)} &= \text{prox}_{\frac{1}{\rho}\mathcal I_{\mathcal C}} (\Delta Y^{(k)}+v^{(k)}) = P_{\mathcal C} (\Delta S^{(k+1)}-v^{(k)}) \notag\\ 
&= \bar E (\bar E^T\bar E)^{-1} \bar E (\Delta S^{(k+1)}-v^{(k)}) = \bar E (\bar E^T\bar E)^{-1} \bar E\Delta S^{(k+1)}\label{eq:ADMMb}
\end{align}
\end{subequations}
\normalsize
\vspace{-5pt}
where  
\small
\begin{equation*}
\begin{split}
\bar f_i(\Delta s^{i}) &= f_i(s^{i,(l)}) + \nabla f_i(s^{i,(l)})^T \Delta s^i + \frac{1}{2}  (\Delta s^i)^T \nabla^2f_i(s^{i,(l)}) \Delta s^i\\
\Delta Y^{k} &= (\Delta y^{1,(k)}, \dots, \Delta y^{N,(k)})\\
v^{k} &= (v^{1,(k)}, \dots, v^{N,(k)})
\end{split}
\end{equation*} 
\normalsize
\vspace*{-10pt}

 and $v^{i,(k)} = v^{i,(k-1)} + (\Delta y^{i,(k)}-\Delta x^{i,(k)})$, which is given by the 5th step of Algorithm \ref{alg:alg3}. Notice that the last equality in~\eqref{eq:ADMMb} holds since if $v^{(0)} = 0$, then $v^{(k)} \in \mathcal N(\bar E (\bar E^T\bar E)^{-1} \bar E)$ for all $k\geq 1$. As we see from~\eqref{eq:ADMMa}, the update for $\Delta S^{(k+1)}$ can be computed in a decentralized manner using $N$ computing agents, where at each iteration every agent $i$ calculates  
\small
\begin{equation}\label{eq:ADMMS}
\begin{split} 
\Delta s^{i,(k+1)} &=  \argmin_{\Delta s^i} \quad \bar f_i(\Delta s^{i}) + \frac{\rho}{2}\| \Delta y^{i,(k)}+v^{i,(k)} - \Delta s^{i}  \|^2\\
& =  \left( \nabla^2f_i(s^{i,(l)}) + \rho I \right)^{-1}\left(\rho(\Delta y^{(k)} + v^{i,(k)}) -  \nabla f_i(s^{i,(l)})\right)
\end{split}
\end{equation}
\normalsize
\vspace*{-12pt}

where $\Delta y^{i,(k)}$ and $v^{i,(k)}$ are locally available to each agent. The update in \eqref{eq:ADMMb} can also be performed distributedly, which is explained as follows.  Let $\Delta z^{(k+1)} = (\bar E^T\bar E)^{-1} \bar E\Delta S^{(k+1)}$. The update rule in~\eqref{eq:ADMMb} can then be rewritten as
\small
\begin{equation}
\begin{split}
\Delta Y^{(k+1)}& = (\Delta y^{1,(k+1)}, \dots, \Delta y^{N,(k+1)})\\ &= \bar E \Delta z^{(k+1)} = (\Delta z_{J_1}^{(k+1)}, \dots, \Delta z_{J_N}^{(k+1)} ).
\end{split}
\end{equation}
\normalsize
\vspace*{-5pt}

Notice that $\bar E^T\bar E = \diag( |\mathcal{I}_1|, \ldots, |\mathcal{I}_{N}|)$, and hence the update for each component, $j$, of $\Delta x$ can be expressed as
\small
\begin{align}\label{eq:ComponentUpdate}
\Delta z^{(k+1)}_j = \frac{1}{|\mathcal{I}_j|}\sum_{q \in \mathcal{I}_j}^{} \left( E_{J_q}^T \Delta s^{q,(k+1)} \right)_j.
\end{align}
\normalsize
\vspace*{-15pt}

As a result, for each agent $i$ to compute $\Delta y^{i,(k+1)} = \Delta z_{J_i}^{(k+1)}$, it needs to communicate to all agents in $\text{Ne}(i) = \left\{ j \ | \ J_i\cap J_j \neq \emptyset  \right\}$ which are referred to as the neighbors of agent $i$. The ADMM-based Newton direction computation can then be summarized as in Algorithm \ref{alg:alg4}.
\vspace{-5pt}
%
%
%
%
\begin{algorithm}[H]
\caption{ADMM-Based Newton Direction Computation}\label{alg:alg4}
\begin{algorithmic}[1]
\small
\State{Given $k= 0$, $\rho > 0$, $\epsilon_{\text{pri}},\epsilon_{\text{dual}}>0$, $\Delta z^{(0)}$ and $v^{(0)}= 0$}
\Repeat
\For{i = 1, 2, \dots,N}
\State {\vspace*{-14pt}\begin{multline*}\Delta s^{i,(k+1)} =  \left( \nabla^2f_i(s^{i,(l)}) + \rho I \right)^{-1}\\ \left( \rho(\Delta z_{J_i}^{(k)} + v^{i,(k)}) - \nabla f_i(s^{i,(l)})\right)\end{multline*}}
\State {Communicate with all agents $r$ belonging to $\text{Ne}(i)$}
\For {all $j \in J_i$}
 \State $\Delta z_j^{(k+1)} = \frac{1}{|\mathcal{I}_j|}\sum_{q \in \mathcal{I}_j}^{} \left( E_{J_q}^T \Delta s^{q,(k+1)} \right)_j$
 \EndFor
 \State $v^{i,(k+1)} = v^{i,(k)} + \left(z_{J_i}^{(k+1)} - s^{i,(k+1)} \right)$
\State  Check whether $\| \Delta z_{J_i}^{(k+1)}-\Delta z_{J_i}^{(k)}\|^2 \leq \epsilon_{\text{dual}}/N$ and $\| \Delta s^{i,(k+1)}-\Delta z_{J_i}^{i,(k+1)}\|^2 \leq \epsilon_{\text{pri}}/N$
\EndFor
\If  {condition in step (10) satisfied for all i = 1, \dots, N} 
\State{$\Delta x_{\nt} =  \Delta z^{(k+1)}$}
\State{$\Delta S_{\nt} = \bar E \Delta z^{(k+1)}$}
\State {Terminate the algorithm} 
\EndIf 
\State{$k = k+1$}
\Until{algorithm is terminated}
\normalsize
 \end{algorithmic}
 \end{algorithm}
\vspace*{-15pt}
Notice that the termination condition of Algorithm \ref{alg:alg4} (based on the 10th and 12th steps of the algorithm) can be established distributedly, provided that all agents declare their status of convergence (step 10). Also observe that the satisfaction of this termination condition implies the satisfaction of the termination condition in Algorithm~\ref{alg:alg3}. There are other ways of establishing convergence of Algorithm \ref{alg:alg4} to a solution (and possibly more efficient, e.g., based on \cite{Iut:12}). However, for the sake of brevity we abstain from discussing such methods in this paper. 

\begin{rem}\label{rem:rem1}
In Algorithm \ref{alg:alg4}, the computational effort at each iteration for each agent is dominated by the update of $\Delta s^{i,(k+1)}$ which requires factorizing the matrix $\nabla^2f_i(s^{i,(l)}) + \rho I$. Notice that in case $\rho$ is chosen to be constant for all $k$, this matrix will also be constant for all $k$. This means that each agent would only need to compute this factorization once at the first iteration of the algorithm and use the precomputed factorization in the remaining iterations. This pre-caching of the factorization significantly reduces the overall computational cost of computing the Newton direction using Algorithm~\ref{alg:alg4}. Notice that even if the penalty parameter is not chosen to be constant, it is still possible to utilize the method presented in \cite[sec. 4.2]{liu:13}, to use pre-caching of factorizations to achieve similar computational efficiency.
\end{rem}

\vspace{-10pt}
\section{Newton's Method with Distributed Step Computation}\label{sec:DNewM}
\vspace{-10pt}
We can now combine Algorithm \ref{alg:alg4} with the Newton method expressed in Algorithm \ref{alg:alg1}, to provide a distributed computational scheme for the problem in \eqref{eq:UnconstGen}. This is presented in Algorithm \ref{alg:alg5}.
\begin{algorithm}[tb]
\caption{ADMM-based Newton Method}\label{alg:alg5}
\begin{algorithmic}[1]
\small
\State{Given $l =0$, $\epsilon_{\nt}>0$, $S^{(0)} \in \mathcal C(\bar E)$, $\Delta x_{\nt}^{(-1)}$, }
\Repeat
\State{Compute $\Delta S_{\nt}^{(l)}$, $\Delta x_{\nt}^{(l)}$ using Algorithm \ref{alg:alg4} with starting point $\Delta z^{(0)} = \Delta x^{(l-1)}_{\nt}$}
\For{i = 1, 2, \dots,N}
\State{Compute $(\lambda^{i,(l)})^2 = (\Delta s_{\nt}^{i,(l)})^T\nabla^2f_i(s^{i,(l)})\Delta s_{\nt}^{i,(l)}$}
\State {Check whether $(\lambda^{i,(l)})^2/2 \leq \epsilon_{\nt}/N$}
\EndFor
\If  {condition in step (6) satisfied for all i = 1, \dots, N} 
\State{$S^\ast = S^{(l)}$}
\State {Terminate the algorithm} 
\EndIf 
\State {Compute the step size $\alpha^{(l)}$ using line search}
\For{i = 1, 2, \dots,N}
\State  $s^{i,(l+1)} =  s^{i,(l)} + \alpha^{(l)}\Delta s^{i,(l)}_{\nt}$
\EndFor
\State{$l = l+1$}
\Until{algorithm is terminated}
\normalsize
\end{algorithmic}
\end{algorithm}
\vspace*{-8pt}
Similar to Algorithm \ref{alg:alg4}, the termination condition for the Newton iterations in Algorithm \ref{alg:alg5} can also be checked in a distributed manner with limited communications among the agents. To this end, each agent needs to compute its local Newton decrement, $\lambda^{i,(l)}$, and declare its local convergence status based on this quantity (condition in step (6)). Although not explained in Algorithm \ref{alg:alg5}, the step size computation in step (11) of the algorithm can also be performed distributedly. For this purpose, each agent would firstly need to perform a back tracking line search based on its local objective function, $\bar f_i(\Delta s^{i})$, and compute a suitable local step size. The step size $\alpha^{(l)}$ will then be chosen as the smallest of the local step sizes, which can be derived using max/min consensus algorithms, \cite{Iut:12}. A similar approach has also been used in \cite{wei:13}. 

\begin{rem} \label{rem:rem2}
Having defined Algorithm \ref{alg:alg5}, some comments are in order. Firstly, notice that in the 3rd step of Algorithm \ref{alg:alg5}, the ADMM iterations for computing the Newton direction, are warm-started using the computed Newton direction from the previous Newton iteration. In case the line search is not done too aggressively this can potentially reduce the number of required ADMM iterations for finding the next Newton direction. Secondly, observe that the computed Newton directions $\Delta S_{\nt}^{(l)} \in \mathcal C(\bar E)$ for all $l \geq 0$ and hence, $S^{(l)} \in \mathcal C(\bar E)$ for all $l \geq 0$. This means that the consistency constraints in \eqref{eq:DEOPb} are always satisfied.   
\end{rem}
So far we have proposed a distributed scheme for solving unconstrained loosely coupled problems. In the upcoming section, we extend the definition of loosely coupled problems to that of constrained ones and show how we can derive similar algorithms to solve such problems.

\vspace*{-10pt}
\section{Constrained Loosely Coupled Optimization Problems}\label{sec:DefCon}
\vspace*{-10pt}
We now extend the definition of loosely coupled problems (provided in Section \ref{sec:Definition}) by first adding convex inequality constraints as
\small
\begin{equation}\label{eq:constGen}
\begin{split}
\minimize & \quad F_1(x) + \dots + F_N(x)\\
\subject & \quad G_i(x) \leq 0  \quad  i = 1, \dots, N
\end{split}
\end{equation}
\normalsize
where we assume that the function pairs $F_i, G_i$ for $i = 1, \dots, N$, are only dependent on a small subset of the elements of the variable $x$. Also let the description of the coupling among the function pairs be described as was for the problem in \eqref{eq:UnconstGen}. Then this problem can be written as
\small
\begin{subequations}\label{eq:ConstDEOP}
\begin{align}
\minimize_{S,x}  & \quad  f_1(s^1) + \dots + f_N(s^N)\\
\subject & \quad   g_i(s^i) \leq 0, \quad i= 1, \dots, N\\
& \quad   \bar E x = S
\end{align}
\end{subequations}
\normalsize
\vspace*{-18pt}

where the functions $f_i$ and $g_i$ are defined similarly as in Section \ref{sec:Definition} for the problem in \eqref{eq:DEOP}. Notice that this problem is in the same format as \eqref{eq:EOPInC} and can be solved using the interior-point method, which accordingly requires solving 
\vspace*{-2pt}
\small
\begin{equation}\label{eq:constGenLogLoose}
\begin{split}
\minimize&  \quad \underbrace{tf_1(s^1) - \log(-g_1(s^1))}_{h_1(s^1)} + \dots +  \underbrace{tf_N(s^N) - \log(-g_N(s^N))}_{h_N(s^N)}\\
\subject & \quad \bar E x = S 
\end{split}
\end{equation}
\normalsize
\vspace*{-15pt}

for a sequence of increasing $t$. For every given $t$ then the problem in \eqref{eq:constGenLogLoose} can be solved using Algorithm \ref{alg:alg5}, with a modification to step $3$ of the algorithm. Particularly, the $4$th step of Algorithm \ref{alg:alg4} (that is used in step 3 of Algorithm~\ref{alg:alg5}) needs to be modified as
\vspace{-7pt}
\small
\begin{multline*}
\Delta s^{i,(k+1)} =  \left( \nabla^2h_i(s^{i,(l)}) + \rho I \right)^{-1}\\ \left( \rho(\Delta z_{J_i}^{(k)} + v^{i,(k)}) - \nabla h_i(s^{i,(l)})\right)
\end{multline*}
\normalsize
\vspace{-2pt}
where 
\vspace*{-5pt}
\small
\begin{multline*}
\nabla^2h_i(s^{i,(l)}) = \nabla^2f_i(s^{i,(l)}) + \frac{1}{g_i(s^{i,(l)})^2} \nabla g_i(s^{i,(l)})\nabla g_i(s^{i,(l)})^T -\\ \frac{1}{g_i(s^{i,(l)})}\nabla^2g_i(s^{i,(l)})  
\end{multline*}
\normalsize
\begin{algorithm}[tb]
\caption{ADMM-based Interior-point Method}\label{alg:alg6}
\begin{algorithmic}[1]
\small
\State{Given $q = 0$, $t^{(0)} >0$, $\mu>1$, $\epsilon_{p}>0$ and feasible $S^{(0)}$}
\Repeat
\State{Compute $S^\ast(t^{q})$ by solving \eqref{eq:EOPInB} using Alg. \ref{alg:alg5} starting at $S^{(q)}$}
\State{$S^{(q+1)} = S^\ast(t^{q})$}
\If  {$m/t^{(q)} < \epsilon_p$} 
\State{$S^\ast = S^{(q+1)}$}
\State {Terminate the iterations} 
\EndIf 
\State{$t^{(q+1)} = \mu t^{(q)}$}
\State{$q = q+1$}
\Until{iterations are terminated}
\normalsize
\end{algorithmic}
\end{algorithm}
and $\nabla h_i(s^{i,(l)}) = \nabla f_i(s^{i,(l)}) - \frac{1}{g_i(s^{i,(l)})} \nabla g_i(s^{i,(l)})$. Notice that after this modification remarks \ref{rem:rem1} and \ref{rem:rem2} still apply. The distributed scheme for solving \eqref{eq:constGen} can then be obtained by combining the modified Algorithm \ref{alg:alg5} with Algorithm \ref{alg:alg2}. This is summarized in Algorithm \ref{alg:alg6}. 
The problem in \eqref{eq:constGen} can be further extended by adding equality constraints. Assuming that the equality constraints also enjoy similar type of coupling as above, the problem can be written as 

\vspace*{-12pt}
\small
\begin{subequations}\label{eq:ConstEDEOP}
\begin{align}
\minimize_{S,x}  & \quad  f_1(s^1) + \dots + f_N(s^N)\\
\subject & \quad   g_i(s^i) \leq 0, \quad i= 1, \dots, N\\
& \quad A^i s^i = b^i \hspace{6mm} i= 1, \dots, N\\
& \quad   \bar E x = S
\end{align}
\end{subequations}
\normalsize
\vspace*{-15pt}

where $A^i \in \mathbb R^{p_i\times |J_i|}$ and $\rank(A^i) = p_i < |J_i|$. This problem can also be solved using Algorithm \ref{alg:alg6} where in its 3rd step, Algorithm \ref{alg:alg5} is applied to 
\small
\begin{equation*}
\begin{split}
\minimize&  \quad \underbrace{tf_1(s^1) - \log(-g_1(s^1))}_{h_1(s^1)} + \dots +  \underbrace{tf_N(s^N) - \log(-g_N(s^N))}_{h_N(s^N)}\\
\subject & \quad A^i s^i = b^i \quad i= 1, \dots, N \\
& \quad \bar E x = S 
\end{split}
\end{equation*}
\normalsize
\vspace*{-15pt}

This in turn requires another modification to Algorithm~\ref{alg:alg4}, that is used in the 3rd step of Algorithm \ref{alg:alg5}. Specifically, the 4th step of Algorithm \ref{alg:alg4} should be changed to 
\small
\begin{multline}\label{eq:ADMMEq}
\Delta s^{i,(k+1)} = \argmin_{A^i \Delta s^i = 0 } \left\{ \nabla h_i(s^{i,(l)})^T \Delta s^i  +  \frac{1}{2}  (\Delta s^i)^T \nabla^2h_i(s^{i,(l)}) \Delta s^i \right. \\ \left.+ \frac{\rho}{2} \left\| z_{J_i}^{(k)} + v^{i,(k)} - \Delta s^i \right\|^2\right\}
\end{multline}
\normalsize
\vspace*{-11pt}

Solving this problem is equivalent to solving the following linear system of equations
\vspace{-5pt}
\small
\begin{multline}
\begin{bmatrix}  \nabla^2h_i(s^{i,(l)}) + \rho I & (A^i)^T \\ A^i & 0 \end{bmatrix}\begin{bmatrix} \Delta s^i \\ u \end{bmatrix}\\ = \begin{bmatrix}\rho(\Delta z_{J_i}^{(k)} + v^{i,(k)})  - \nabla h_i(s^{i,(l)}) \\ 0 \end{bmatrix}.
\end{multline}
\normalsize
\vspace{-11pt}

The computational cost for this is dominated by the cost for calculating the factorization of the coefficient matrix. Notice that, similar to \eqref{eq:ADMMS}, this coefficient matrix is also constant within the ADMM iterations and hence, the comments made in Remark \ref{rem:rem1} still apply. Even if $\rho$ is varying similar techniques as in \cite{liu:13} can be applied.
\begin{rem}
Notice that the resulting solution from Algorithm \ref{alg:alg6}, $S^\ast$, satisfies the local constraints, i.e., $g_i(s^{i,\ast}) \leq 0$ and $A^is^{i,\ast} = b^i$ for all $i= 1, \dots, N$. However, $S^\ast$ does not necessarily satisfy the consistency constraints. In fact it is possible to upper-bound the consistency error in $S^\ast$ as follows. Assume that $S^{(0)}$ is both locally feasible and consistent. Then the consistency error of $S^{(1)}$ can be computed as
\small
\begin{equation}
\begin{split}
e_c(1) &= \left\| S^{(1)} - \bar E (\bar E^T \bar E)^{-1}\bar E^T S^{(1)}\right\|^2 \\&= \sum_{j=0}^{l_\text{max}(1)} (\alpha^{(j,1)})^2\left\| \Delta S^{j,1} - \bar E (\bar E^T \bar E)^{-1}\bar E^T\Delta S^{(j,1)} \right\|^2 \\
& \leq \sum_{j=0}^{l_\text{max}(1)} (\alpha^{(j,1)})^2 \epsilon_{\text{pri}},
\end{split}
\end{equation}
\normalsize
\vspace*{-10pt}

where $\Delta S^{(j,1)} $ and $\alpha^{(j,1)}$ denote the Newton direction the step size used in the $j$th Newton iteration in the 1st iteration of the interior-point method, respectively. Also $l_\text{max}(1)$ denotes the number of required Newton iterations to converge to $S^{(1)}$. Now assuming that Algorithm \ref{alg:alg6} requires $q_{\text{max}}$ iterations to converge, we can bound the consistency error of $S^{\ast}$ as $e_c(q_{\text{max}}) \leq \sum_{q=1}^{q_{\text{max}}} \sum_{j=0}^{l_\text{max}(q)} (\alpha^{(j,q)})^2 \epsilon_{\text{pri}}$.
\end{rem}
\begin{rem}
The proposed methods in this paper are closely related to inexact Newton methods, where the Newton direction is computed in a distributed manner. This indicates that the methods described in this paper, can potentially suffer from slow convergence when the iterates are very close to the optimal solution. The similarities between the two classes of methods can also give us insight on the convergence properties of the proposed methods and how to avoid slow convergence near optimal solution.
\end{rem}
\vspace*{-10pt}
\section{Conclusions}\label{sec:Con}
\vspace*{-8pt}

In this paper, we proposed distributed optimization algorithms for solving unconstrained and constrained loosely coupled optimization problems. These algorithms are based on the Newton and interior-point methods where we used the inherent structure in the problem to distribute the required computations within these methods. Particularly, by exploiting the coupling in the problem and using proximal splitting methods we showed how the Newton direction can be computed in a distributed manner. 
\newline
Notice that the proposed methods in this paper, all require a feasible starting point. As a future research directions, we intend to extend the proposed algorithms to also accommodate infeasible starting points. Also we intend to investigate primal-dual interior-point methods.

\vspace{-7pt}
\small
\bibliography{IEEETrans}

\end{document}